\documentclass[12pt]{article}
\usepackage{amsmath,amssymb,amsbsy,amsfonts,amsthm,latexsym,amsopn,amstext,
amsopn,amstext,amsxtra,euscript,amscd}
\begin{document}
\newfont{\teneufm}{eufm10}
\newfont{\seveneufm}{eufm7}
\newfont{\fiveeufm}{eufm5}
%
%
\newfam\eufmfam
                             \textfont\eufmfam=\teneufm
\scriptfont\eufmfam=\seveneufm
                             \scriptscriptfont\eufmfam=\fiveeufm

%
%
\def\frak#1{{\fam\eufmfam\relax#1}}
%


\def\bbbr{{\rm I\!R}} 
\def\bbbc{{\rm I\!C}} 
\def\bbbm{{\rm I\!M}}
\def\bbbn{{\rm I\!N}} 
\def\bbbf{{\rm I\!F}}
\def\bbbh{{\rm I\!H}}
\def\bbbk{{\rm I\!K}}
\def\bbbl{{\rm I\!L}}
\def\bbbp{{\rm I\!P}}
\newcommand{\lcm}{{\rm lcm}}
\def\bbbone{{\mathchoice {\rm 1\mskip-4mu l} {\rm 1\mskip-4mu l}
{\rm 1\mskip-4.5mu l} {\rm 1\mskip-5mu l}}}
\def\bbbc{{\mathchoice {\setbox0=\hbox{$\displaystyle\rm C$}\hbox{\hbox
to0pt{\kern0.4\wd0\vrule height0.9\ht0\hss}\box0}}
{\setbox0=\hbox{$\textstyle\rm C$}\hbox{\hbox
to0pt{\kern0.4\wd0\vrule height0.9\ht0\hss}\box0}}
{\setbox0=\hbox{$\scriptstyle\rm C$}\hbox{\hbox
to0pt{\kern0.4\wd0\vrule height0.9\ht0\hss}\box0}}
{\setbox0=\hbox{$\scriptscriptstyle\rm C$}\hbox{\hbox
to0pt{\kern0.4\wd0\vrule height0.9\ht0\hss}\box0}}}}
\def\bbbq{{\mathchoice {\setbox0=\hbox{$\displaystyle\rm
Q$}\hbox{\raise
0.15\ht0\hbox to0pt{\kern0.4\wd0\vrule height0.8\ht0\hss}\box0}}
{\setbox0=\hbox{$\textstyle\rm Q$}\hbox{\raise
0.15\ht0\hbox to0pt{\kern0.4\wd0\vrule height0.8\ht0\hss}\box0}}
{\setbox0=\hbox{$\scriptstyle\rm Q$}\hbox{\raise
0.15\ht0\hbox to0pt{\kern0.4\wd0\vrule height0.7\ht0\hss}\box0}}
{\setbox0=\hbox{$\scriptscriptstyle\rm Q$}\hbox{\raise
0.15\ht0\hbox to0pt{\kern0.4\wd0\vrule height0.7\ht0\hss}\box0}}}}
\def\bbbt{{\mathchoice {\setbox0=\hbox{$\displaystyle\rm
T$}\hbox{\hbox to0pt{\kern0.3\wd0\vrule height0.9\ht0\hss}\box0}}
{\setbox0=\hbox{$\textstyle\rm T$}\hbox{\hbox
to0pt{\kern0.3\wd0\vrule height0.9\ht0\hss}\box0}}
{\setbox0=\hbox{$\scriptstyle\rm T$}\hbox{\hbox
to0pt{\kern0.3\wd0\vrule height0.9\ht0\hss}\box0}}
{\setbox0=\hbox{$\scriptscriptstyle\rm T$}\hbox{\hbox
to0pt{\kern0.3\wd0\vrule height0.9\ht0\hss}\box0}}}}
\def\bbbs{{\mathchoice
{\setbox0=\hbox{$\displaystyle     \rm S$}\hbox{\raise0.5\ht0\hbox
to0pt{\kern0.35\wd0\vrule height0.45\ht0\hss}\hbox
to0pt{\kern0.55\wd0\vrule height0.5\ht0\hss}\box0}}
{\setbox0=\hbox{$\textstyle        \rm S$}\hbox{\raise0.5\ht0\hbox
to0pt{\kern0.35\wd0\vrule height0.45\ht0\hss}\hbox
to0pt{\kern0.55\wd0\vrule height0.5\ht0\hss}\box0}}
{\setbox0=\hbox{$\scriptstyle      \rm S$}\hbox{\raise0.5\ht0\hbox
to0pt{\kern0.35\wd0\vrule height0.45\ht0\hss}\raise0.05\ht0\hbox
to0pt{\kern0.5\wd0\vrule height0.45\ht0\hss}\box0}}
{\setbox0=\hbox{$\scriptscriptstyle\rm S$}\hbox{\raise0.5\ht0\hbox
to0pt{\kern0.4\wd0\vrule height0.45\ht0\hss}\raise0.05\ht0\hbox
to0pt{\kern0.55\wd0\vrule height0.45\ht0\hss}\box0}}}}
\def\bbbz{{\mathchoice {\hbox{$\sf\textstyle Z\kern-0.4em Z$}}
{\hbox{$\sf\textstyle Z\kern-0.4em Z$}}
{\hbox{$\sf\scriptstyle Z\kern-0.3em Z$}}
{\hbox{$\sf\scriptscriptstyle Z\kern-0.2em Z$}}}}
\def\ts{\thinspace}

\newtheorem{theorem}{Theorem}
\newtheorem{lemma}[theorem]{Lemma}
\newtheorem{claim}[theorem]{Claim}
\newtheorem{cor}[theorem]{Corollary}
\newtheorem{prop}[theorem]{Proposition}
\newtheorem{definition}{Definition}
\newtheorem{question}[theorem]{Open Question}

\def\squareforqed{\hbox{\rlap{$\sqcap$}$\sqcup$}}
\def\qed{\ifmmode\squareforqed\else{\unskip\nobreak\hfil
\penalty50\hskip1em\null\nobreak\hfil\squareforqed
\parfillskip=0pt\finalhyphendemerits=0\endgraf}\fi}

\def\cA{{\mathcal A}}
\def\cB{{\mathcal B}}
\def\cC{{\mathcal C}}
\def\cD{{\mathcal D}}
\def\cE{{\mathcal E}}
\def\cF{{\mathcal F}}
\def\cG{{\mathcal G}}
\def\cH{{\mathcal H}}
\def\cI{{\mathcal I}}
\def\cJ{{\mathcal J}}
\def\cK{{\mathcal K}}
\def\cL{{\mathcal L}}
\def\cM{{\mathcal M}}
\def\cN{{\mathcal N}}
\def\cO{{\mathcal O}}
\def\cP{{\mathcal P}}
\def\cQ{{\mathcal Q}}
\def\cR{{\mathcal R}}
\def\cS{{\mathcal S}}
\def\cT{{\mathcal T}}
\def\cU{{\mathcal U}}
\def\cV{{\mathcal V}}
\def\cW{{\mathcal W}}
\def\cX{{\mathcal X}}
\def\cY{{\mathcal Y}}
\def\cZ{{\mathcal Z}}

\newcommand{\comm}[1]{\marginpar{%
\vskip-\baselineskip 
\raggedright\footnotesize
\itshape\hrule\smallskip#1\par\smallskip\hrule}}





\def\ve{\varepsilon}

\hyphenation{re-pub-lished}

\def\ord{{\mathrm{ord}}}
\def\Nm{{\mathrm{Nm}}}
\renewcommand{\vec}[1]{\mathbf{#1}}

\def \F{{\bbbf}}
\def \L{{\bbbl}}
\def \K{{\bbbk}}
\def \Z{{\bbbz}}
\def \N{{\bbbn}}
\def \Q{{\bbbq}}
\def\E{{\mathbf E}}
\def\bH{{\mathbf H}}
\def\G{{\mathcal G}}
\def\O{{\mathcal O}}
\def\cS{{\mathcal S}}
\def \R{{\bbbr}}
\def\Fp{\F_p}
\def \fp{\Fp^*}
\def\\{\cr}
\def\({\left(}
\def\){\right)}
\def\fl#1{\left\lfloor#1\right\rfloor}
\def\rf#1{\left\lceil#1\right\rceil}

\def\Zm{\Z_m}
\def\Zt{\Z_t}
\def\Zp{\Z_p}
\def\Um{\cU_m}
\def\Ut{\cU_t}
\def\Up{\cU_p}

\def\ep{{\mathbf{e}}_p}

\def\e{{\mathbf{e}}}

\def \Prob{{\mathrm {}}}

\def\LC{{\cL}_{C,\cF}(Q)}
\def\LCn{{\cL}_{C,\cF}(nG)}
\def\Mrs{\cM_{r,s}\(\F_p\)}

\def\Fbar{\overline{\F}_q}
\def\Fn{\F_{q^n}}
\def\En{\E(\Fn)}

\def \hatf{\widehat{f}}

\def\mand{\qquad \mbox{and} \qquad}

\def\MOV{{\bf{MOV}}}

\title{\bf On the Distribution of Pseudopowers}

\author{
{\sc Sergei V. Konyagin} \\
             {Department of Mechanics and Mathematics}\\
{Moscow State University} \\
{Moscow, 119992, Russia} \\
{\tt konyagin@ok.ru}\\
             \and
{\sc Carl Pomerance}\\
{Department of Mathematics}\\ {Dartmouth College}\\
{Hanover, NH 03755-3551, USA} \\
{\tt carlp@gauss.dartmouth.edu} \\
\and
{\sc Igor E. Shparlinski}\\
             {Department of Computing}\\
{Macquarie University}\\
{ Sydney, NSW 2109,
Australia}\\
{\tt igor@ics.mq.edu.au}
}

\date{}

\maketitle

\begin{abstract}
An $x$-pseudopower to base $g$ is a positive integer which is not a power
of $g$ yet is so modulo $p$ for all primes $p\le x$.
We improve an upper bound for the least such number due to
E.~Bach, R.~Lukes, J.~Shallit, and H.~C.~Williams.
The method is based on a combination of
some bounds of exponential sums with new results
about the average behaviour of the multiplicative order
of  $g$ modulo prime numbers.
\end{abstract}


\section{Introduction}

Let $g$ be a fixed integer with $|g|\ge 2$.
Following E.~Bach, R.~Lukes, J.~Shallit, and
H.~C.~Williams~\cite{BLSW}, we say that  an
integer $n>0$ is  an {\it $x$-pseudopower to  base $g$\/}
if $n$ is not a power of
$g$ over the integers but is
a power of $g$  modulo all primes $p\le x$,
that is, if for all primes $p\le x$ there exists an integer $e_p\ge
0$ such that
$n\equiv g^{e_p} \pmod p $.

Denote by $q_g(x)$ the least $x$-pseudopower to base $g$.

A well-known result of A.~Schinzel~\cite{Schin}
asserts that if  $f$ and $g>0$ are integers, such that
$f\neq g^k$ for all integers $k\ge 0$, then for infinitely many primes $p$
the congruence $g^x\equiv f \pmod p$ does not have solutions in
nonnegative integers $x$. Therefore,
$$
q_g(x)\to\infty, \qquad x\to\infty.
$$
E.~Bach, R.~Lukes, J.~Shallit and H.~C.~Williams~\cite{BLSW}
have shown that if the Riemann hypothesis holds for Dedekind zeta
functions, then
there is a constant $A>0$, depending only on $g$, such that
$$
q_g(x)\ge\exp(A\sqrt x/(\log x)^2).
$$

On the other hand, 
if
$$
M_x  = \prod_{p\le x} p
$$
is the product of all primes
$p\le x$, then $q_g(x)\le 2M_x+1$ when $x\ge 2$.  Indeed,
both $M_x+1$ and $2M_x+1$ are $\equiv g^0\pmod p$ for all
primes $p\le x$ and evidently not both can be powers of $g$.
The prime number theorem implies that $M_x=e^{(1+o(1))x}$, so we have
\begin{equation}
\label{eq:CRT Bound}
q_g(x)\le e^{(1+o(1))x}, \qquad x\to\infty.
\end{equation}
Though the inequality $q_g(x)\le2M_x+1$ cannot be improved in general (consider
the case $g=M_x+1$), if $g$ is fixed or $|g|$ is not too large
compared with $x$, there is a chance to improve the bound~\eqref{eq:CRT Bound}.
Supported by numerical data, a heuristic argument is given in~\cite{BLSW}
suggesting that $q_g(x)$ for fixed $g$
is about $\exp(c_gx/\log x)$, where $c_g>0$.
We obtain a more modest upper bound valid for $|g|\le x$
as well as several more results about the distribution of
$x$-pseudopowers to base $g$.

For an integer $m$ we use $\Z_m$ to denote the residue
ring modulo $m$. Now,  for a prime $p$,  we denote by
$\cU_{g,p}$ the subset of $\Z_p$ generated by powers of $g$ modulo $p$,
that is
$$
\cU_{g,p}=\{n\in \Z_p\ :\   n\equiv g^k \pmod p \ \ \text{for some
nonnegative}\ k  \in \Z\}.
$$
Clearly, if $\gcd(g,p)=1$ then $\cU_{g,p}$ is a subgroup
of $\Z_p^*$, while if $p\mid g$, then $\cU_{g,p}=\{0,1\}$.

We consider the set
$$
\cW_g(x) = \{n\in[0,M_x)\ :\  n \in \cU_{g,p} \
\text{for all primes}\ p \le x\}.
$$
The set $\cW_g(x)$ consists of both the $x$-pseudopowers to base $g$
that lie below $M_x$ and the true powers of $g$ in this range.  (In the
case that $M_x\mid g$, the set $\cW_g(x)$ also contains 0, but
we shall be assuming that $|g|\le x$ and $x$ is large, so that this case 
does not occur.)  
The number of true powers of $g$ below $M_x$ is $O(x)$, which 
turns out to be minuscule in comparison to $\#\cW_g(x)$.

We first get a good lower bound for $\#\cW_g(x)$.  Then
we estimate exponential sums with elements of $\cW_g(x)$ and
use these bounds to derive some uniformity-of-distribution results
for elements of $\cW_g(x)$.  Our estimate for $q_g(x)$ follows
from these results.

\section{Our approach and results}

Our approach is based on a combination of two
techniques:

\begin{itemize}
\item recent bounds of exponential sums over reasonably
small subgroups of the  multiplicative group $\Z_p^*$
due to Heath-Brown and Konyagin~\cite{HBK};
\item Lower bounds on multiplicative orders on average which
we derive from upper bounds of R.~C.~Baker
and G.~Harman~\cite{BaHa1,BaHa2} (which are summarised in~\cite{Harm}) for the
Brun--Titchmarsh inequality on average.
\end{itemize}

We do not try to obtain numerically the best results,
rather we concentrate on the exposition of our main
ideas. Certainly with more work and numerical calculations
one can get more precise results. Furthermore,  any further
advance in our knowledge on the above two topics would
immediately lead to further progress on this problem as well.

For prime $p\nmid g$, let $l_g(p)=\#\cU_{g,p}$, the multiplicative order of
$g$ modulo $p$.
We also put $l_g(p) =1$ for $g \equiv 0 \pmod p$.
We now define the product
\begin{equation}
\label{eq:Rproduct}
R_g(x)
    = \prod_{p\le x}  l_g(p).
\end{equation}
The Chinese remainder theorem implies that
$$
\#\cW_g(x)=\prod_{p\le x}\#\cU_{g,p}.
$$
Further, for $p\mid g$, we have 
$\#\cU_{g,p}=2=2l_g(p)$.  Thus, if $\gcd(g,M_x)$ has exactly
$k$ prime factors,
\begin{equation}
\label{eq:R and W}
\#\cW_g(x)=2^kR_g(x)\ge R_g(x).
\end{equation}

Note that $R_g(x)^{1/\pi(x)}$ is the geometric mean of $l_g(p)$ for $p\le x$
and so has some independent interest.  Our first result gives a lower bound
for $R_g(x)$ and so, via \eqref{eq:R and W}, gives a lower bound for $\#\cW_g(x)$.

\begin{theorem}
\label{thm:Rlowerbound}
For $x$ sufficiently large and for
$g$ an integer with $2\le|g|\le x$, we have
$$
\#\cW_g(x)\ge R_g(x)\ge \exp(\eta x)
$$
where $\eta=0.58045$.
\end{theorem}

We put $\e(u)=\exp(2\pi iu)$ and define exponential sums
$$
S_{a,g}(x) = \sum_{n\in \cW_g(x) }\e(an/M_x) .
$$
\begin{theorem}
\label{thm:exp sum}
For $x$ sufficiently large and for any integers $a,g$ with 
$2\le|g|\le x$, we have
$$
|S_{a,g}(x)|\le \#\cW_g(x)\gcd(a,M_x)\exp(-\gamma x)
$$
where
$$\gamma=0.11286.$$
\end{theorem}

For a positive  integer   $h\le M_x$, let
$N_g(x,h)$ denote the number of members of $\cW_g(x)$ below $h$.
Using some standard arguments, we derive from
our estimates of the sums $S_{a,g}(x)$:

\begin{theorem}
\label{thm:unif distr}
For $x$ sufficiently large, we have for any integers $g$ and
$h$ with $2\le |g|\le x$ and $1\le h\le M_x$,
$$N_g(x,h)= \#\cW_g(x) \frac{h}{M_x}+E_g(x,h)$$
where
$$
|E_g(x,h)|\le\#\cW_g(x)\exp(-\gamma x)$$
and where $\gamma$ is as in Theorem~\ref{thm:exp sum}.
\end{theorem}

In particular, we improve~\eqref{eq:CRT Bound} to
$$
q_g(x)\le e^{0.88715x}
$$
for $x$ sufficiently large and $|g|\le x$.  Indeed, if we take $h=e^{0.88715x}$
in Theorem~\ref{thm:unif distr}, then that result implies that there
are at least $\frac12\#\cW_g(x)h/M_x$ numbers in $\cW_g(x)$ below $h$.
Together with Theorem~\ref{thm:Rlowerbound} this implies that there
are more than $e^{.4675x}$ members of $\cW_g(x)$ below $h$.
But there are only $O(x)$ numbers below $h$ that are true powers of $g$,
so there are many members of $\cW_g(x)$ below $h$ that are $x$-pseudopowers
to base $g$.

\section{Proof of Theorem~\ref{thm:Rlowerbound}}
\label{sec:gmmo}

It is well known, see~\cite{ErdMur,Ford,IndlTim,Papp},
that $l_g(p)\ge x^{1/2}$ for all but
$o(x/\log x)$ primes $p\le x$. Thus for $R_g(x)$,
given by~\eqref{eq:Rproduct},
we immediately obtain
\begin{equation}
\label{eq:Rtriv}
R_g(x) \ge \exp(x/2 + o(x)).
\end{equation}
We now obtain a more accurate estimate for $R_g(x)$.

Let $P(m)$ denote the largest prime divisor of $m\ge 2$
(with the convention $P(1) = 0$).
We use $\pi(x,y)$ to denote the number of
primes $p\le x$ with $P(p-1)\le y$ and define the constant
\begin{equation}
\label{eq:const c}
c =  \liminf \pi(x,x^{1/2})/\pi(x).
\end{equation}

\begin{lemma}
\label{lem:R and c}
For the product $R_g(x)$, given by~\eqref{eq:Rproduct}, we have
$$
R_g(x) \ge \exp\( \frac{1 + c}{2}  x + o(x)\),
$$
where $c$ is given by~\eqref{eq:const c}.
\end{lemma}

\begin{proof}
Let $\cP_0$ be the set of primes $p\le x$ with $l_g(p)\le x^{1/2}$,
let $\cP_1$ be the set of primes $p\le x$ with $l_g(p)>x^{1/2}$
and $P(p-1)> x^{1/2}$, and let $\cP_2$ be the set of all other primes $p\le x$.

We simply ignore the contribution from primes in $\cP_0$ (which,
as we have mentioned, is $\exp(o(x))$ anyway).

For each $p\in \cP_1$,  since $l_g(p)\mid p-1$ and $(p-1)/P(p-1)<x^{1/2}$,
we have $P(p-1)\mid l_g(p)$.  Thus,
\begin{equation}
\label{P1sum-1}
\sum_{p\in\cP_1}\log l_g(p)\ge \sum_{p\in\cP_1}\log P(p-1)=
\sum_{x^{1/2}<q\le x}\pi(x;q,1)\log q~+o(x),
\end{equation}
where $q$ runs over primes and $\pi(x;k,b)$ denotes the number of primes
$p\le x$ with
$p \equiv b \pmod k$.  Indeed, each $q$ in the indicated
range corresponds to $\pi(x;q,1)$ primes $p\le x$ with $P(p-1)=q$,
and almost all primes $p$ so counted in the sum are in $\cP_1$.
It follows  from
the Bombieri--Vinogradov theorem and the Brun--Titchmarsh inequality
(see~\cite[Theorems~6.6 and 17.1]{IwKow}) that
$$
\sum_{q\le x^{1/2}}\pi(x;q,1)\log q= (1/2+o(1))x,
$$
and since $\sum_{q\le x}\pi(x;q,1)\log q=(1+o(1))x$, we have
\begin{equation}
\label{Gold}
\sum_{x^{1/2}<q\le x}\pi(x;q,1)\log q= (1/2+o(1))x,
\end{equation}
as noted by M.~Goldfeld~\cite{Gold}.  We thus have
from \eqref{P1sum-1} that
\begin{equation}
\label{P1sum-2}
\sum_{p\in\cP_1}\log l_g(p)\ge (1/2+o(1))x.
\end{equation}

We now consider the contribution from primes in $\cP_2$.  For each
such prime $p$
we have $l_g(p)\ge x^{1/2}$, so that
\begin{equation}
\label{P2sum}
\sum_{p\in\cP_2}\log l_g(p)\ge \frac12\log x\sum_{p\in\cP_2}1
=\frac12\pi(x,x^{1/2})\log x +o(x).
\end{equation}

The bounds~\eqref{P1sum-2} and~\eqref{P2sum}, together with~\eqref{eq:const c},
imply that
$$
\sum_{p\le x}\log l_g(p)  \ge (1/2 + c/2 + o(1)) x,
$$
which concludes the proof.
\end{proof}

There is probably little doubt that
$$
c = \rho(2) = 1 - \log 2  = 0.3068\ldots\,,
$$
where $\rho(u)$ is the Dickman--de~Bruijn function (see~\cite{Ten}),
however proving this seems to be inaccessible
by present methods; see~\cite{BFPS,Pom,PomShp} where more general
conjectures about $\pi(x,y)$ are discussed.
Note that in~\cite{Pom} we have the inequality
$$
\pi(x,x^{1/2})\ge(1-4\log(5/4)+o(1))x/\log x,
$$
so that $c \ge 0.107425\ldots $.
The key tool in~\cite{Pom} is a result of C.~Hooley~\cite{Hool} from 1973.
Using more modern tools we now obtain a  larger value of $c$.

For $1/2\le u<1$ let
$C(u)$ denote a monotone nondecreasing function
such that for any $\ve>0$ and $A>0$, we have
\begin{equation}
\label{BTavg}
\pi(x;k,b)\le (C(u)+\ve)\frac{x}{\varphi(k)\log x}
\end{equation}
for all integers $k\le x^u$ but for at most $x^u/\log^Ax$ exceptions,
for all $b$ coprime to $k$ for allowable values of $k$, and
for all $x\ge x_0(A,\ve)$.
H.~L.~Montgomery and R.~C.~Vaughan~\cite{MoVa} have a
version of the Brun--Titchmarsh theorem which allows one to take
$C(u)=2/(1-u)$ with no exceptional values of $k$ and with $\ve=0$,
see also~\cite[Theorem~8.1]{Harm} or~\cite[Section~6.8]{IwKow}).
But allowing a small exceptional set
as indicated here then permits one to get smaller values of $C(u)$.  This
is the arena of ``the Brun--Titchmarsh theorem on average."
The key results we use are due to {\'E}.~Fouvry~\cite{Fouv} and
R.~C.~Baker and G.~Harman~\cite{BaHa1,BaHa2}.
(There are many other contributors to this subject,
we refer to~\cite{Harm} for more details and further references).

For a monotone nondecreasing function $C(u)$ satisfying \eqref{BTavg},
let us define $\vartheta_C$ by the  equation
\begin{equation}
\label{thetaprop}
\int_{1/2}^{\vartheta_C} C(u)\,du = 1/2.
\end{equation}
(Note that for any monotone nondecreasing function $C(u)$ the
integral is well defined.)

We now use the approach of~\cite{Pom} to show the following lower
bound on $c$.

\begin{lemma}
\label{lem:c and C(u)}
For the constant $c$ given by~\eqref{eq:const c}, we
have
$$
c \ge 1 - \int_{1/2}^{\vartheta_C}\frac{C(u)}{u}\,du
$$
where $C(u)$ is an arbitrary
      monotone nondecreasing function satisfying~\eqref{BTavg}
and $\vartheta_C$ is defined by~\eqref{thetaprop}.
\end{lemma}

\begin{proof} Let
$$
H(x,t)=\sum_{x^{1/2}<q\le t}\pi(x;q,1)\log q
$$
where $q$ runs over primes.
Thus, by partial summation, we have
\begin{equation}
\label{nonsmooths}
\pi(x)-\pi(x,x^{1/2})=\sum_{x^{1/2}<q\le x}\pi(x;q,1)
       =\frac{H(x,x)}{\log x}+\int_{x^{1/2}}^x\frac{H(x,t)}{t\log^2t}\,dt.
\end{equation}
Using~\eqref{Gold}, the first term on the right in~\eqref{nonsmooths}
is $(1/2+o(1))x/\log x$, so it remains
to get a good upper bound for the integral.

Using the inequality~\eqref{BTavg}, partial summation, and the prime number
theorem, we have
\begin{equation}
\label{eq:Hxtineq}
H(x,t)\le x\int_{1/2}^{\log t/\log x}C(u)\,du+o(x).
\end{equation}
Thus, for any value of $\vartheta\in(1/2,1)$ we have
\begin{equation}
\label{truncint}
\int_{x^{1/2}}^{x^\vartheta}\frac{H(x,t)}{t\log^2t}\,dt
\le x\int_{x^{1/2}}^{x^\vartheta}\frac{1}{t\log^2t}\int_{1/2}^{\log
t/\log x}C(u)\,du\,dt
+o(x/\log x).
\end{equation}
By a change of variables and an interchange of the order of
integration, the double
integral is equal to
\begin{eqnarray*}
\int_{x^{1/2}}^{x^\vartheta}\frac{1}{t\log^2t}\int_{1/2}^{\log t/\log
x}C(u)\,du\,dt
& = &  \int_{1/2}^{\vartheta}C(u)
\int_{x^u}^{x^{\vartheta}}\frac{1}{t\log^2t} dt\,du\\
& = &  \int_{1/2}^{\vartheta}C(u)  \int_{u \log x}^{\vartheta \log x
}\frac{1}{v^2} dv\,du\\
& = &\frac1{\log x}\int_{1/2}^\vartheta C(u)\(\frac1u-\frac1\vartheta\)\,du.
\end{eqnarray*}

Thus, from~\eqref{truncint} we have
$$
\int_{x^{1/2}}^{x^\vartheta}\frac{H(x,t)}{t\log^2t}\,dt
\le \frac{x}{\log x}\int_{1/2}^\vartheta
C(u)\(\frac1u-\frac1\vartheta\)\,du+o(x/\log x).
$$
Using $H(x,t)\le H(x,x)=(1/2+o(1))x$ (see~\eqref{Gold}), we then have
for any $\vartheta\in(1/2,1)$
that
\begin{equation*}
\begin{split}
\int_{x^{1/2}}^x\frac{H(x,t)}{t\log^2t}\,dt
&=\int_{x^{1/2}}^{x^\vartheta}\frac{H(x,t)}{t\log^2t}\,dt
+\int_{x^\vartheta}^x\frac{H(x,t)}{t\log^2t}\,dt\\
\le\frac{x}{\log x}&\int_{1/2}^\vartheta C(u)\(\frac1u-\frac1\vartheta\)\,du
+\frac{x}{2\log x}\(\frac1\vartheta-1+o(1)\)
\end{split}
\end{equation*}
which we rewrite as
\begin{equation}
\begin{split}
\label{Hintineq}
\int_{x^{1/2}}^x\frac{H(x,t)}{t\log^2t}& \,dt\\
       \le\frac{x}{\log x}&\( \int_{1/2}^\vartheta  \frac{C(u)}{u}\,du
- \frac{1}{\vartheta }  \int_{1/2}^\vartheta C(u) du + \frac{1}{2\vartheta }
- \frac{1}{2} +o(1) \).
\end{split}
\end{equation}
If we choose  $\vartheta = \vartheta_C$ defined by~\eqref{thetaprop},
then  using~\eqref{nonsmooths} and~\eqref{Hintineq},  we obtain
$$
\pi(x)-\pi(x,x^{1/2})\le\frac{x}{\log
x}\int_{1/2}^{\vartheta_C}\frac{C(u)}{u}\,du+o(x/\log x)
$$
which concludes the proof.
\end{proof}

We now use known results on the possible choices of the
function $C(u)$ in~\eqref{BTavg}, as summarised in~\cite{Harm},
to obtain a lower bound for $c$.

\begin{lemma}
\label{lem:Const c}
For the constant $c$ given by~\eqref{eq:const c}, we have
$$
c >0.160901.
$$
\end{lemma}

\begin{proof}
For $u \in [0.51, 0.56]$, we define $C(u)$ as
a step-wise monotonically nondecreasing function
whose values at $u = 0.533$ and  $u = 0.5 + 0.005j$, $j = 1, \ldots, 12$
are   given by G.~Harman in~\cite[Theorem~8.2]{Harm} as $C(0.533) = 2$
and in~\cite[Table~8.1]{Harm} as:
\begin{center}
\begin{tabular}{|l|l|  |l|l| |l| l| |l|l||l|l| }
\hline
      $u$ & $C(u)$& $u$ & $C(u) $ & $u$ & $C(u)$ & $u$ & $C(u) $ & $u$
& $C(u) $\\
\hline
0.515 & 1.223 &  0.525 & 1.75& 0.535 & 2.09 & 0.545 & 2.47  & 0.555 & 2.76\\
0.52 & 1.632 & 0.53 &1.82 & 0.540 &2.25 & 0.55 & 2.66 & 0.56 & 2.88 \\
\hline
\end{tabular}
\end{center}

For other values of $u$,
we also use  analytic expressions which are  due to
R.~C.~Baker and G.~Harman~\cite{BaHa1,BaHa2} and
{\'E}.~Fouvry~\cite{Fouv}. These results are also
presented in~\cite[page~184]{Harm} (for  $u \in [0.5, 0.51)$)
and in~\cite[Theorem~8.4]{Harm} ($u \in [17/32, 5/7]$):

\begin{itemize}

\item for $0.5 \le u< 0.51$, we have $C(u) = 1+ 150(u-1/2)^2$;

\item for $17/32 < u\le 4/7$, we have $C(u) =
14/(12-13u)-\log(4(1-u)/3u)$ (in fact we use it only
for $0.56< u \le 4/7$);

\item for $4/7< u\le 3/5$ we have $C(u) =14/(12-13u)$;

\item for $3/5 < u\le 5/7$ we have $C(u)  = 8/(3-u)$.
\end{itemize}

With this we compute (using {\sl Mathematica})
\begin{eqnarray*}
\int_{0.5}^{0.51} C(u)\,du=0.01005\,, & &
\int_{5.1}^{0.56} C(u)\,du  =  0.107405\,, \\
\int_{0.56}^{4/7} C(u) \,du
\approx 0.034177\,,  & &
\int_{4/7}^{3/5} C(u) \,du
\approx 0.091260\,, \\
\int_{3/5}^{.6759} C(u)\,du
\approx 0.257087\,, & &
\int_{3/5}^{.67591} C(u)\,du
\approx 0.257121\,,
\end{eqnarray*}
where the approximations are rounded to 6 decimal places.

Therefore
$$
\int_{1/2}^{0.6759} C(u)\,du
<0.49999,
\qquad \int_{1/2}^{0.67591}C(u)\,du
>0.50001,
$$
and we see that for our choice of $C(u)$,
$$
0.6759<\vartheta_C < 0.67591.
$$

We also compute
\begin{eqnarray*}
\int_{0.5}^{0.51} \frac{C(u)}{u}\,du<0.019902\,, & &
\int_{0.51}^{0.56} \frac{C(u)}{u}\,du<0.199610\,,\\
\int_{0.56}^{4/7} \frac{C(u)}{u}\,du <  0.060412\,, & &
\int_{4/7}^{3/5} \frac{C(u)}{u}\,du < 0.155787\,, \\
\int_{3/5}^{0.67591} \frac{C(u)}{u}\,du <  0.403388\,.
\end{eqnarray*}
Therefore,
$$
\int_{1/2}^{\vartheta_C}\frac{C(u)}{u}\,du <
\int_{1/2}^{0.67591}\frac{C(u)}{u}\,du  <  0.839099\,,
$$
so that with Lemma \ref{lem:c and C(u)} the result follows.
\end{proof}

Theorem \ref{thm:Rlowerbound} now follows from Lemmas
\ref{lem:R and c} and \ref{lem:Const c}.

\section{Proof of Theorem \ref{thm:exp sum}}
\label{section exp sum}

Let
\begin{equation}
\label{eq:Tproduct}
T_g(x) = \prod_{\substack{p\le x\\
l_g(p) > p^{2/3}} }(l_g(p)p^{-2/3}).
\end{equation}
We now show that the exponential sum $S_{a,g}(x)$ is
related to $T_g(x)$ and the product $R_g(x)$ defined
in \eqref{eq:Rproduct}.

\begin{lemma}
\label{lem:prelim}  For any integer $a$, we have
$$
|S_{a,g}(x)| \le
\#\cW_g(x) d \exp\(x/4 + o(x) \) R_g(x)^{-5/8} T_g(x)^{-3/8}
      $$
as $x\to\infty$, where $d=\gcd(a,M_x)$.
\end{lemma}

\begin{proof}
By the Chinese remainder theorem we see that
\begin{equation}
\label{eq:CRT}
S_{a,g}(x)=\sum_{n\in \cW_g(x) }\e(an/M_x)  =\prod_{p\le x}
\sum_{n\in \cU_{g,p}}\e(a_pn/p)
\end{equation}
where $a_p\in\Zp$ is determined by the condition
$$a_p(M_x/p)\equiv a \pmod {M_x}.$$

If $p\nmid ag$, the bound of D.~R.~Heath-Brown and S.~V.~Konyagin~\cite{HBK}
applies which gives the estimate
\begin{equation}
\label{eq:HBK}
\left| \sum_{n\in \cU_{g,p}}\e(a_pn/p) \right|\le C l_g(p)^{3/8}
p^{1/4}.
\end{equation}
for some absolute constant $C > 1$.
We also recall  the well-known bound
\begin{equation}
\label{eq:Kor}
\left| \sum_{n\in \cU_{g,p}}\e(a_pn/p) \right| \le    p^{1/2}
\end{equation}
(provided $p\nmid ag$),
which is better than~\eqref{eq:HBK} for $l_g(p) > p^{2/3}$,
see~\cite[Theorem~3.4]{KoSh}.

For the set  $\cP_0$ of primes $p \le x$
with $p\mid ag$ we estimate the exponential sums over $\cU_{g,p}$
trivially as $2l_g(p)$.

For the set  $\cP_1$ of primes
with  $p\nmid ag$ and $l_g(p) \le
p^{2/3}$ we use the bound~\eqref{eq:HBK}.

Finally, for the set  $\cP_2$ of primes
with  $p\nmid ag$ and $l_g(p) >
p^{2/3}$ we use the bound~\eqref{eq:Kor}.

Thus, substituting these bounds in~\eqref{eq:CRT}, we obtain
\begin{equation}
\label{Sagest}
|S_{a,g}(x)|\le 2^{\#\cP_0}C^{\#\cP_1}\prod_{p\in\cP_0 }l_g(p)
\prod_{p\in\cP_1}l_g(p)^{3/8}p^{1/4} \prod_{p\in\cP_2}p^{1/2}.
\end{equation}
We majorize the first two factors in~\eqref{Sagest} as
$e^{O(\pi(x))}=e^{o(x)}$.
The first product in \eqref{Sagest} may be restricted to the primes
$p\le x$ which
divide $a$, and since $l_g(p)<p$, this product is bounded by $\gcd(a,M_x)=d$.
Let $\cQ_1,\cQ_2$ be the same as $\cP_1,\cP_2$ but without the
restriction that $p\nmid ag$.  Thus, the three products in
\eqref{Sagest} are at most
\begin{eqnarray*}
d\prod_{p\in\cQ_1}l_g(p)^{3/8}p^{1/4}\prod_{p\in\cQ_2}p^{1/2}
&=&d\prod_{p\le x}l_g(p)^{3/8}p^{1/4}\prod_{p\in\cQ_2}l_g(p)^{-3/8}p^{1/4}\\
&=&dR_g(x)^{3/8}M_x^{1/4}T_g(x)^{-3/8}.
\end{eqnarray*}
Thus, the result follows from~\eqref{Sagest}, the prime number theorem
in the form $M_x=e^{(1+o(1))x}$, and the inequality \eqref{eq:R and W}.
\end{proof}

Using the elementary bound \eqref{eq:Rtriv} together with Lemma~\ref{lem:prelim}
and the trivial bound $T_g(x)\ge 1$
already gives a nontrivial estimate on the sums $S_{a,g}(x)$, namely
$$
|S_{a,g}(x)| \le \#\cW_g(x) \gcd(a,M_x)\exp\(-x/16 + o(x) \) .
$$
Using Theorem \ref{thm:Rlowerbound} in place of \eqref{eq:Rtriv}
and still using only $T_g(x)\ge 1$ we get
$$
|S_{a,g}(x)| \le \#\cW_g(x) \gcd(a,M_x)\exp\(-0.11278x \)
$$
for all large $x$.
We now obtain a nontrivial estimate for $T_g(x)$,
which in turn implies a slightly better estimate for $S_{a,g}(x)$.

\begin{lemma}
\label{lem:prod T}
For the product $T_g(x)$ given by~\eqref{eq:Tproduct},
a function $C(u)$ satisfying~\eqref{BTavg}, and $\vartheta_C>\frac{2}{3}$
defined by~\eqref{thetaprop} we have
$$
T_g(x) \ge \exp\(x\int_{2/3}^{\vartheta_C}\left(1-\frac2{3u}\right)C(u)\,du
+o(x)\).
     $$
\end{lemma}

\begin{proof}
Let $\cP$ be the set of primes $p\le x$ with $l_g(p)>x^{1/2}$ and
$P(p-1)> x^{2/3}$.  Similarly to the proof of Lemma~\ref{lem:R and c} 
we have $l_g(p)>p^{2/3}$ for all $p\in\cP$ and so
\begin{equation}
\begin{split}
\label{P1sum-3}
\log T_g(x)&\ge\sum_{p\in\cP}(\log l_g(p)-\frac{2}{3}\log p)\\
&\ge\sum_{x^{2/3}<q\le x}\pi(x;q,1)(\log q-\frac{2}{3}\log x)+o(x),
\end{split}
\end{equation}
where $q$ runs over primes.

Next, we follow the proof of Lemma~\ref{lem:c and C(u)}.
By partial summation, we have
\begin{equation}
\label{nonsmooths2}
\begin{split}
\sum_{x^{2/3}<q\le x}\pi(x;q,1)&(\log q-\frac{2}{3}\log x)\\
&=\frac13H(x,x)-\frac{2 \log x}{3}\int_{x^{2/3}}^x\frac{H(x,t)}{t\log^2t}\,dt.
\end{split}
\end{equation}
Using~\eqref{eq:Hxtineq} as in the argument for~\eqref{Hintineq}, and
recalling \eqref{thetaprop}, we get
\begin{equation*}
\begin{split}
\int_{x^{2/3}}^x\frac{H(x,t)}{t\log^2t}  \,dt ~\le ~&\frac{x}{\log x}
\Biggl( \int_{2/3}^{\vartheta_C}    \frac{C(u)}{u}\,du
    +\frac32\int_{1/2}^{2/3} C(u) du \\
& \qquad\qquad \qquad-~
\frac{1}{\vartheta_C }  \int_{1/2}^{\vartheta_C} C(u) du +
\frac{1}{2\vartheta_C} - \frac{1}{2} +o(1)\Biggr) \\
~=~& \frac{x}{\log
x}\(\int_{2/3}^{\vartheta_C}\frac{C(u)}{u}\,du + \frac{3}{2}
\int_{1/2}^{2/3} C(u) du-\frac12+o(1)\).
\end{split}
\end{equation*}
Combining this with~\eqref{P1sum-3} and~\eqref{nonsmooths2},
and then using~\eqref{Gold},
     we complete the proof.
\end{proof}

Using the estimates for the function $C(u)$ as discussed in the
proof of Lemma~\ref{lem:Const c} we can now get an explicit estimate
for $T_g(x)$.
\begin{lemma}
\label{lem:Test}
For the product $T_g(x)$, given by \eqref{eq:Tproduct}, and $x$ sufficiently
large, we have
$$
T_g(x)\ge\exp(0.000217x).
$$
\end{lemma}

\begin{proof}
This follows immediately from Lemma~\ref{lem:prod T}, the
estimate $\vartheta_C>0.6759$ seen in the proof of Lemma~\ref{lem:Const c},
and the formula $C(u)=8/(3-u)$ for the range $[3/5,5/7]$
also seen in the proof of Lemma~\ref{lem:Const c}.
\end{proof}

We now have Theorem~\ref{thm:exp sum} by using,
in the inequality of Lemma~\ref{lem:prelim},
our estimate for $R_g(x)$ in Theorem~\ref{thm:Rlowerbound}
and our estimate for $T_g(x)$ in Lemma~\ref{lem:Test}.

\section{Proof of Theorem~\ref{thm:unif distr}}

Using that for any integer $m \ge 1$ we have
$$
\sum_{a=0}^{m-1} \e (au/m) =
\left\{ \begin{array}{ll}
0,& \quad \mbox{if}\ u\not \equiv 0 \pmod m, \\
m,& \quad \mbox{if}\ u \equiv 0 \pmod m,
\end{array} \right.
$$
(which follows from the formula for the sum of
a geometric progression) we write
$$ N_g(x,h) = \sum_{n \in  \cW_g(x)}\sum_{k=0}^{h-1}
\frac{1}{M_x}\sum_{a=0}^{M_x-1} \e \(a(n-k)/M_x\).
    $$
Changing the order of summation and separating the term
$\#\cW_g(x)h/M_x$ corresponding to $a=0$ we derive
\begin{equation}
\label{eq:N Delta}
\left|N_g(x,h) - \#\cW_g(x)\frac{h}{M_x} \right| \le \frac{1}{M_x} \Delta
\end{equation}
where
$$
\Delta = \sum_{a=1}^{M_x-1} |S_{a,g}(x)|  \left| \sum_{k=0}^{h-1}
\e\(ak/M_x\)\right|.
    $$
For each $d \mid M_x$ with $d < M_x$ we now collect together the terms
with $\gcd(a,M_x)=d$ and also apply Lemma~\ref{lem:prelim}, getting the
estimate
\begin{eqnarray*}
\Delta & \le &\#\cW_g(x)    \exp\(x/4 + o(x) \) R_g(x)^{-5/8} T_g(x)^{-3/8}
    \\
& &\qquad \qquad \qquad \qquad\qquad
\sum_{\substack{d < M_x\\ d \mid M_x}} d \sum_{\substack{a =1\\
\gcd(a,M_x)=d}}^{M_x-1}
\left|\sum_{k=0}^{h-1} \e\(ak/M_x\)\right|\\
& \le &\#\cW_g(x)    \exp\(x/4 + o(x) \) R_g(x)^{-5/8} T_g(x)^{-3/8}\\
& &\qquad \qquad \qquad \qquad\qquad
\sum_{\substack{d < M_x\\ d \mid M_x}} d \sum_{b =1}^{M_x/d -1}
\left|\sum_{k=0}^{h-1} \e\(bk/(M_x/d)\)\right|.
\end{eqnarray*}
We now recall that for any integers $m \ge 2$ and $1 \le b < m$,
    we have the bound
$$
\left|\sum_{k=0}^{h-1} \e\(bk/m\)\right| \ll \frac{m}{\min\{b, m-b\}}
$$
which again follows from the
formula for the sum of
a geometric progression, see~\cite[Bound~(8.6)]{IwKow}.
This implies that
$$
\sum_{b=1}^{m-1} \left|\sum_{k=0}^{h-1} \e\(bk/m\)\right|
\ll m \log m.
$$
Thus
$$
\Delta  \le \#\cW_g(x) M_x \exp\(x/4 + o(x) \) R_g(x)^{-5/8} T_g(x)^{-3/8}
$$
where we used that
$$
\sum_{\substack{d < M_x\\ d \mid M_x}} 1\le 2^{\pi(x)}
= \exp\( o(x) \).
$$
Substituting this bound in~\eqref{eq:N Delta},
we obtain
$$
\left|N_g(x,h) - \#\cW_g(x)\frac{h}{M_x} \right| \le
\#\cW_g(x)
\exp\(x/4 + o(x) \) R_g(x)^{-5/8} T_g(x)^{-3/8} .
$$
Theorem~\ref{thm:unif distr} now follows from our estimates
for $R_g(x)$ and $T_g(x)$ in Sections~\ref{sec:gmmo} and
\ref{section exp sum}, respectively.

\section{Remarks}

Using better estimates for
$C(u)$ that already exist, it is possible to get a larger
value of $\vartheta_C$ and consequently better numbers in
Lemmas~\ref{lem:Const c} and~\ref{lem:Test}.  
In particular in~\cite{BaHa2} and~\cite{Harm} a method of computing a somewhat smaller
function $C$ satisfying~\eqref{BTavg} is described leading to 
$\vartheta_C>0.677$.
Using this value of $\vartheta_C$ in our estimate for
$T_g(x)$ allows us to replace $0.000217$ with $0.000272$.
The changes in the estimate for $c$ in Lemma~\ref{lem:Const c}
depend much more intrinsically on the better estimates for $C(u)$
that support a value of $\vartheta_C$ that is greater than $0.677$;
we have not worked this out.

Certainly if more information about the possible choice of the
function $C(u)$ becomes available, one can immediately obtain
even better numerical estimates for the constant $c$ and thus
improve the results of
Theorems~\ref{thm:Rlowerbound}, \ref{thm:exp sum}, and~\ref{thm:unif distr}.

Another avenue for improvement could come with our estimate
for $l_g(p)$ when $P(p-1)\le\sqrt{x}$.  We used the estimate
$l_g(p)\ge\sqrt{x}$ for almost all such
primes $p\le x$.  It follows from~\cite[Theorem~6]{Ford} of K.~Ford
that there is some $\varepsilon>0$ such that for a positive
proportion of these primes we have $l_g(p)\ge x^{1/2+\varepsilon}$.
Having a version of this theorem with explicit constants would
allow a numerical improvement in our Lemma~\ref{lem:R and c}
and thus an improvement in our principal results.

It is very plausible that the technique of~\cite[Chapter~7]{KoSh}
can be used to improve our bound on $q_g(x)$ (but not the
bounds of Theorems~\ref{thm:exp sum} and~\ref{thm:unif distr}).
However adjusting this technique to the case of composite
moduli and then tuning it to accomodate in an optimal way
our current knowledge of the behaviour of $l_g(p)$ may take
significant efforts. 

Finally,  we recall that under the Generalised Riemann
Hypothesis we have $l_g(p) = p^{1 + o(1)}$ for almost all
primes $p$, see~\cite{ErdMur,KP,Papp}, which immediately gives
$$
R_g(x) = \exp(x + o(x))  \mand  T_g(x) = \exp(x/3 + o(x)) .
$$
In turn, this means that one can take any $\gamma < 1/2$ in
Theorems~\ref{thm:exp sum} and~\ref{thm:unif distr} and one
has $q_g(x)\le e^{x/2+o(x)}$.

\section*{Acknowledgments}  We thank G. Harman for sending
us a pre-publication file for Chapter 8 of his book~\cite{Harm}.
S.V.K.\ gratefully acknowledges support from Grants 05-01-00066
from the Russian Foundation for Basic Research and NSh-5813.2006.1.
C.P.\ gratefully acknowledges support from NSF grants DMS-0401422 and DMS-0703850.
I.E.S.\ gratefully acknowledges support from ARC grant DP0556431.

\end{document}